\documentclass[reqno]{amsart}

\usepackage{amssymb,amscd}
\usepackage{color}
\usepackage{graphicx}

\theoremstyle{plain}
\newtheorem{theorem}{Theorem}[section]
\newtheorem{lemma}{Lemma}[section]

\newtheorem*{theorem*}{Theorem}

\theoremstyle{definition}
\newtheorem{definition}{Definition}[section]

\theoremstyle{remark}
\newtheorem*{remark}{Remark}

\newcommand{\N}{\ensuremath{\mathbb{N}}}
\newcommand{\Z}{\ensuremath{\mathbb{Z}}}

\newcommand{\R}{\ensuremath{\mathbb{R}}}

\newcommand{\KK}{\ensuremath{\mathcal{K}}}
\newcommand{\LL}{\ensuremath{\mathcal{L}}}
\newcommand{\QQ}{\ensuremath{\mathcal{Q}}}
\newcommand{\Sph}{\ensuremath{\mathcal{S}}}
\newcommand{\UU}{\ensuremath{\mathcal{U}}}

\newcommand{\Wo}{\ensuremath{\mathcal{W}_{O}}}
\newcommand{\Wf}{\ensuremath{\mathcal{W}_{F}}}
\newcommand{\Wr}{\ensuremath{\mathcal{W}_{R}}}

\newcommand{\Gfa}{\ensuremath{\hat{f}_{b}}}
\newcommand{\GL}{\ensuremath{GL_d(\R)}}
\newcommand{\card}{\ensuremath{\#}}

\newcommand{\za}{\ensuremath{\alpha}}

\newcommand{\zd}{\ensuremath{\delta}}

\newcommand{\zve}{\ensuremath{\varepsilon}}
\newcommand{\zvf}{\ensuremath{\varphi}}
\newcommand{\zC}{\ensuremath{\Psi}}
\newcommand{\zg}{\ensuremath{\gamma}}
\newcommand{\zG}{\ensuremath{\Gamma}}
\newcommand{\zo}{\ensuremath{\omega}}
\newcommand{\zl}{\ensuremath{\lambda}}

\newcommand{\zO}{\ensuremath{\Omega}}

\newcommand{\zr}{\ensuremath{\rho}}

\newcommand{\zp}{\ensuremath{\psi}}

\newcommand{\beq}{\begin {equation}}
\newcommand{\eeq}{\end{equation}}

\newcommand{\benu}{\begin {enumerate}}
\newcommand{\eenu}{\end{enumerate}}

\newcommand{\bite}{\begin {itemize}}
\newcommand{\eite}{\end{itemize}}

\newcommand{\bdef}{\begin {definition}}

\newcommand{\diam}{\ensuremath{\text{diam}}}
\newcommand{\supp}{\ensuremath{\text{supp}}}
\newcommand{\dist}{\ensuremath{\text{dist}}}

\begin{document}


\title
{Density of the set of generators of wavelet systems}


\author[C.Cabrelli]{Carlos~Cabrelli}
\address{Departamento de
Matem\'atica \\ Facultad de Ciencias Exactas y Naturales\\ Universidad
de Buenos Aires\\ Ciudad Universitaria, Pabell\'on I\\ 1428 Capital
Federal\\ ARGENTINA\\ and CONICET, Argentina}
\email[Carlos~Cabrelli]{cabrelli@dm.uba.ar} 

\author[U.Molter]{Ursula~M.~Molter}
\email[Ursula~M.~Molter]{umolter@dm.uba.ar}

\begin{abstract}
Given a function $\psi$ in $ \LL^2(\R^d)$, the affine  (wavelet) system generated by $\psi$, associated
to an invertible matrix $a$ and a lattice $\zG$, is the collection of functions
$\{|\det a|^{j/2} \psi(a^jx-\gamma): j \in \Z, \gamma \in \zG\}$.
In this article we prove that the set of functions generating affine systems that are a Riesz basis of $ \LL^2(\R^d)$ is dense
in $ \LL^2(\R^d)$.

We also prove that a stronger result is true for affine systems that are a frame of $ \LL^2(\R^d)$.
In this case we show that the generators associated to a fixed but arbitrary dilation are a dense set.

Furthermore, we analyze the orthogonal case in which we prove that the set of generators of orthogonal (not necessarily complete) affine systems, that are compactly supported in frequency,
 are dense in the unit sphere of $ \LL^2(\R^d)$ with the induced metric. As a byproduct we
introduce the $p$-Grammian of a function and prove a convergence result of this Grammian as a function of the lattice.
This result gives  insight in the problem of oversampling of affine systems.
\end{abstract}

\keywords{Wavelet Set, affine systems, Riesz basis wavelets, Wavelet Frames}
\subjclass{Primary:42C40, Secondary: 42C30}

\date{\today}

\maketitle


\section{Introduction}
 Let $\psi$ be an $\LL^2(\R^d)$-function, $a$ an invertible $d\times d$ matrix  and $\zG$ a lattice in $\R^d$.
The {\em wavelet system} or {\em affine system}  $(\psi, a, \zG)$,  generated by $\psi$ and associated to $a$ and $\zG$, is the 
collection of functions,

\beq
(\psi, a, \zG) = \left\{D_{a^j}T_\zg \psi: j\in \Z, \zg \in \zG\right\}.
\eeq
Here  $T_y$ and $D_a$ denote the the unitary operators in  $\LL^2(\R^d)$
defined by
$(T_y \psi)(x) = \psi(x-y), y \in \R^d$ and $
(D_a \psi)(x) = |\det a|^{1/2} \psi(ax), a \in \GL.$

Affine systems have been studied in depth during the last 25 years mainly because of their importance in
applications. In addition they proved to be very useful  in a variety  of theoretical
problems.
On the other hand they were studied in the context of Hilbert spaces where the translation and dilation
operators were replaced by a general group of unitary operators \cite{GLT93, DL98}.
Generalizations of wavelet systems in  $\LL^2(\R^d)$ with translations not necessarily on a lattice
and using different dilations were also considered \cite{ACM04}, \cite{GLLWW04}.

One of the relevant questions  about these systems is whether the collection $(\psi, a, \zG)$
form an orthonormal basis, a Riesz basis or a frame of the space $\LL^2(\R^d)$.

The construction of affine systems with specific prescribed properties is a difficult problem. Usually it is accomplished  imposing conditions on the generators. This has been one of the core problem from the beginning of  wavelet research \cite{Dau88, Mey88, Mey92, BF94, CHM04} (and references therein).
Still there are many unanswered questions and open problems in the study of affine systems.

One way to get a better understanding of  these systems is considering the set of generators
of affine systems having some particular structure, and trying to answer global questions
about these sets.
For example, wavelet systems forming a tight frame of the space have been completely characterized by
a set of equations imposing conditions on  the Fourier transform of the generators 
\cite{CCMW02, CS00, Bow03}.
Another example is  the problem of connectivity of the set of generators which has received considerable attention lately, \cite{Wut98, Spe99}.

In this article, motivated by a question posed by  David Larson \cite{Lar05}, we study the problem of density
of the set of generators of affine systems.
We prove (Theorem~\ref{wav-set}) that the set  of generators forming Riesz bases is dense in  $\LL^2(\R^d)$
when we allow diagonal matrices and arbitrary lattices. This theorem gives an even stronger result, since the generators in the dense set  have Fourier transform supported on {\em  wavelet-sets} 
and as a consequence have orthogonal dilations. 
The proof is obtained as a combination of the techniques used in the proof of density for wavelet frames (section~\ref{frames}) and
the theory of wavelets sets that has been developed recently by different groups of researches
\cite{BL01, BMM99, BS02, BS04, DLS97, DLS98, ILP98, SW98, Wan02, Zak96}.

The  question about density of generators of wavelet frames, is also answered positively. We obtain a very general result,
that confirms the flexibility in their construction.
If a dilation matrix $a$ is chosen arbitrarily, the set of frame
generators associated to the dilation $a$ is dense in $\LL^2(\R^d)$ (Theorem~\ref{teo-2}).
The main tool here is the use of the general scheme in frame construction that appears in 
\cite{ACM04}.

Finally in section~\ref{ort-sys} we study the case of orthogonal affine systems.
Since generators of orthonormal systems  have norm one the question here is 
whether they are dense in the  unit  sphere of $\LL^2(\R^d)$.
It is easy to see that this is not the case when either the dilation or the lattice is fixed
(see section~\ref{ort-sys}).

We first show that for a given function of norm one, we can always find a lattice $\zG$
and a function with orthonormal translates in that lattice, that is close to the original
function.
We then prove that a dilation can be chosen  in such a way that the associated affine system
is orthonormal. 
This proves the density of generators of orthonormal wavelet systems 
(not necessarily complete in $\LL^2(\R^d)$) (Theorem~\ref{teo-osys}).

The strategy here  is the study of the behavior of the Grammian of the generator, as a function of the lattice.
We obtain an interesting result on the convergence of the p-Grammian. We prove that
the $p$-Grammian of a function converges point-wise
and in $\LL^p([0,1]^d)$ to the constant function $1$ when the lattice expands. 
This result is a statement about oversampling. 

However the question whether  complete orthonormal affine systems are dense
in the sphere remains unanswered.

\section{Notation}

We will denote by $\mu$ the Lebesgue measure in $\R^d$ and by $GL_d(\R)$ 
the usual group of invertible matrices in $\R^{d\times d}$. By a {\em lattice} $\zG$ we
mean $\zG = c\Z^d$, where $c\in \GL$.

The Fourier transform of a function $f \in \LL^2(\R^d)$ is
\beq \notag
\hat{f}(\zo) = \int_{\R^d} f(x)e^{-2\pi i x \zo} dx.
\eeq
For a measurable $\zO \subset \R^d$,  $K_{\zO}$ will denote the functions in $\LL^2(\R^d)$ with support in $\zO$,
\beq \label{k-omega}
K_{\zO} = \{ f \in \LL^2(\R^d): \supp(f) \subseteq \zO\}.
\eeq
We will say that a matrix $a$ is {\em expansive} if $a \in \GL$ and $|\zl| > 1$ for all  eigenvalues 
$\zl$ of $a$.

Since affine systems depend on the matrix and the lattice involved we consider the following
sets of generators: 

\beq \notag
\Wo(a,\zG)  := \left\{\psi \in \LL^2(\R^d): (\psi, a, \zG)\
\text{is an {\em orthonormal basis (onb)} of}\ \LL^2(\R^d)\right\}
\eeq
and

\begin{align*}
\Wo(a) &:= \bigcup_{c\in \GL} \Wo(a, c\Z^d) \\
\Wo(\zG) &:= \bigcup_{a\in \GL} \Wo(a,\zG) \\
\Wo &:= \bigcup_{a, c \in \GL} \Wo(a, c\Z^d).
\end{align*}

Similarly, we use the notations $\Wr$ or $\Wf$ for the corresponding generators of wavelet Riesz basis and frames respectively.

\section{Density of the set of frame generators}

\label{frames}
In this section we will prove that for a fixed expansive matrix $a$, the set $\Wf(a)$ is dense in $ \LL^2(\R^d)$. 

In what follows $\diam(S)$  will denote the diameter of a set $S \subset \R^d$ and
$B(0,t)$ the ball in $\R^d$ centered at $0$ and with radius $t$.

Let  $X$ be a discrete set in $\R^d$. The {\em gap} $\rho$ of $X$  is defined as:
$$\rho=\rho(X)=\sup_{x\in \R^d} {\inf_{\zg\in
X}|x-\zg|}.
$$
The set $X$ is said to be {\em separated} if
$\inf_{\zg \ne \zg'} |\zg - \zg' |>0.$

For a separated set $X$ Beurling \cite{Beu66} proved the following result:
\begin {theorem*}[Beurling]
Let $X\subset \R^d$ be separated, and $\Omega=B(0,r)$. If
$r\rho<1/4$, then $\{e^{-2\pi i \zg \zo}\chi_{\Omega}(\zo): \zg \in X\}$
is a frame for $\mathcal {K}_\Omega$.
\end {theorem*}
For a very clear exposition of some of the Beurling density results see
\cite{BW99}.

We need now the following result which is a particular case of the more general theorems in \cite{ACM04}. We include its proof here for completeness of this presentation.

\begin{theorem}
\label {ACM}
Let $a$ be an expansive matrix and $h$ a bounded compactly supported function such that $|h(\omega)| > c > 0$, for almost 
every $\omega \in \UU :=  \supp(h)$, and $0 \not\in \supp(h)$.
If $X \subset \R^d$ is a separated set such that the gap 
$\rho(X) < \frac{1}{4R}$ where $R$ is such that $U \subset B(0,R)$
then  
$$
\left\{D_{a^j}T_{\zg} \psi: j\in \Z, \zg \in X \right\}
$$
 is a frame of $ \LL^2(\R^d)$ where $\hat \psi = h$.
\end{theorem}
\begin{proof}
By Beurling's theorem we know that 
$$ \{e^{-2\pi i \zg \zo}\chi_{B(0,R)}(\zo): \zg \in X\}$$
is a frame of $K_{B(0,R)}$ with some frame bounds $0 < m \leq M < +\infty$. In particular 
$$ \{e^{-2\pi i \zg \zo}\chi_{\UU}(\zo): \zg \in X\}$$
is a frame of $K_{\UU}$. Using that $D_a$ is a unitary operator and calling $b=(a^{-1})^t$, we conclude that for each $j \in \Z$
$$ \{|\det(b)|^{j/2}e^{-2\pi i b^j\zg \zo}\chi_{\UU}(b^j\zo): \zg \in X\}$$
is a frame of $K_{b^{-j}\UU}$ with the same bounds $m, M$. 

Further, by Proposition 5.9 in \cite{ACM04} there exist $0<p\leq P<+\infty$ such that
$$
p < \sum_{j \in \Z} |h(b^{j}\zo)|^2 < P, \ \text{a.e.} \ \zo \in \R^d.
$$
Now, since $\supp(h) = \UU$, if we set $h_j(\zo) = h(b^j\zo)$ then $\supp(h_j) = b^{-j}\UU$.

Given now $f \in \LL^2(\R^d)$  and calling $f_j = \overline{h_j} f$ we will see that
\begin{equation}\label{fj}
p\|f\|^2  \leq \  \sum_j \|f_j\|^2 \leq P\|f\|^2,
\end{equation}
for
\begin{align*}
p\|f\|^2 & = \ \int p|f|^2 \  \leq   \int \sum_j |h_j|^2 |f|^2  =
\ \sum_j \int |h_j f|^2= \sum_j \| f_j\|^2,
\end{align*}
and the upper inequality can be obtained similarly.
Using that $f_j \in K_{b^{-j}\UU}$ we have
\beq \label{*}
m \|f_j\|^2     \leq
       \ \sum_{\zg} |<f_j(\zo),|\det(b)|^{j/2}e^{-2\pi i b^j \zg \zo}\chi_{b^{-j}U}(\zo)>|^2
       \leq
\ M\|f_j\|^2.
\eeq
In addition, since $\supp(f_j) = b^{-j}\UU$, we have that
\beq \label{***}
 <f_j(\zo),|\det(b)|^{j/2}e^{-2\pi i b^j \zg \zo}\chi_{b^{-j}\UU}(\zo)> = 
<f(\zo),|\det(b)|^{j/2}e^{-2\pi i b^j \zg \zo}h_j(\zo)>. \eeq
Using \eqref{***} in \eqref{*} and summing in $j$ by \eqref{fj} we get
\begin{align*}
pm\|f\|^2 & \leq 
\ \sum_j\sum_{zg} |<f(\zo),|\det(b)|^{j/2}e^{-2\pi i b^j \zg \zo}h_j(\zo)>|^2 \ \leq \ PM\|f\|^2.
\end{align*}
The claim is now a consequence of Plancherel's Theorem.
\end{proof}
Now we are ready to state and proof the main result of this section.

\begin{theorem} \label{teo-2}
The set $\Wf(a)$ is dense in $ \LL^2(\R^d)$, for every
 $d\times d$ expansive matrix $a$.
 \end{theorem}

\begin{proof}
Assume that a function $f \in  \LL^2(\R^d)$ and a positive number $\varepsilon$ are given.
We want to approximate $f$  by a wavelet frame function $\psi$.
We will do this by constructing  the  Fourier transform of $\psi.$
\begin{itemize}

\item Select $g \in  \LL^2(\R^d)$ such that $\hat g$ is a continuous function and $\|\hat f - \hat g\|_2  < \zve/2$.

\item Choose $0< r < R < \infty$ such that

$$
aB(0,r)  \subset B(0,R) \ \ 
\text{and } \ \ 
\int_{\R^d \setminus \UU} |\hat g(\zo)|^2 d\zo < \frac{\varepsilon^2}{8},
$$
where $\UU= \{\omega: r \leq |\omega| \leq R\}$.

\item   Choose $\lambda > 0$, such that ${\displaystyle 4 \zl^2 \mu(\UU) <\frac{\zve^2}{8}}$,
and let 
\beq \label{elambda}
E_{\lambda} :=\{\zo \in \R^d: |\hat g(\zo)| > \lambda \}.
\eeq
\item Define a function $h$ in the following way:

\begin{equation} \label{h-def}
h(\omega) := \begin{cases}
{\hat g}(\omega) & \omega \in \UU\cap E_{\lambda}\\
\lambda & \omega \in \UU \setminus E_{\lambda}\\
0 & \text{else}
\end{cases}.
\end{equation}
\item 
Choose $X$ to be any separated set with gap $\zr(X) < \frac{1}{4R}$. (Note that $X$ can be chosen to be a lattice.)
\end{itemize}

We will now see that we are indeed under the hypothesis of Theorem~\ref{ACM}.
Since 
\beq \notag
 \QQ :=  aB(0,r) \setminus B(0,r) \subset B(0,R) \setminus B(0,r) = \UU
\eeq
we note that by Lemma 5.11 of \cite{ACM04} $\{a^j \QQ\}$ is a covering of $\R^d\setminus \{0\}$
and therefore $\{a^j U\}$ is also a covering of $\R^d\setminus \{0\}$.

Finally, $h$ and $X$ were defined to match the hypothesis of  Theorem~\ref{ACM}.
Therefore we can apply that Theorem to conclude that $(\psi, a, X)$ is a frame for $\LL^2(\R^d)$, with
$\psi = \hat h$.

Furthermore, 
\begin{align}
\|\hat g - h \|_2^2 & = \int_{\UU} |\hat g - h|^2  + \int_{\R^d \setminus \UU} |\hat g - h|^2 
 \notag \\
& < \int_{\UU\setminus E_\lambda} |\hat g - h|^2 d\zo + \frac{\zve^2}{8} \notag \\
& \leq 4\zl^2 \mu(\UU\setminus E_\lambda) + \frac{\zve^2}{8} < \frac{\zve^2}{4}
\notag
\end{align}
where the last inequality comes from the choice of $\zl$.

Finally
\beq \notag
\|f - \psi\|_2 \leq \|\hat f - \hat g\|_2 + \|\hat g - h\|_2 < \zve.
\eeq
\end{proof}

\section{The Riesz basis case}

In this section we will prove that the set of generators for Riesz basis is dense
when we allow diagonal matrices and arbitrary lattices. Since the requirement of building a Riesz basis is much stronger than to form a frame, the construction of the function requires more subtle techniques than for the previous case.

The idea is again to approximate the given function on the Fourier side 
with an appropriate function. We will use the 
construction for the frame wavelets and adapt it to the Riesz basis case. The main difficulty here lies in the fact that the elements of the affine system need to be independent, which forces the construction to be more involved.

It is worth to remark here, that the construction in the previous section was very general and {\em stable} in the sense that we had freedom to move things a little bit - and still obtain a satisfactory result. In this section, the construction is tight and very adapted to the function we want to approximate.

The method we use, is based on a standard {\em wavelet-set} construction. Wavelet sets have been studied and developed by many groups (for references see the Introduction).
We adapt the construction of \cite{Zak96} and \cite{BS02} to obtain functions supported on a wavelet set, but which are not constant. In this way we obtain a Riesz basis, instead of a wavelet basis.

In what follows we will say that two sets $A$ and $B$ in $\R^d$ are {\em almost disjoint} if $\mu(A\cap B) = 0$. 
\begin{lemma}\label{riesz}
Let $\Omega \subset \R^d$ be a set of finite measure.
If $\{\lambda_k\}_{k \in \Z} \subset \R^d$ satisfies that 
$\left\{e^{2\pi i \lambda_k\zo} \chi_{\zO}(\zo) : k \in \Z\right\}$
is a Riesz basis for $\KK_{\zO}$ with bounds $A$ and $B$ and $h \in \LL^2(\R^d)$ 
satisfies that 
$
0<p<|h(\zo)|<P<+\infty
$
then 
\beq\notag \{h(\zo) e^{2\pi i \lambda_k\zo}: k \in \Z\}\eeq
is a Riesz basis for $\KK_{\Omega}$, with Riesz bounds $pA\mu(\zO)$ and $PB\mu(\zO)$.

 If $a \in \GL$ and $\zO$ satisfies that $\cup_{j \in \Z} a^j \zO = \R^d$ up to a set of zero measure, with the union being almost disjoint, and
$\{g_k: k \in \Z\}$ is a Riesz basis for $\KK_{\Omega}$, then
\beq
\{D_a^j g_k: k,j \in \Z\}
\notag \eeq
is a Riesz basis for $\LL^2(\R^d)$ with the same bounds.
\end{lemma}
\begin{proof}
The first assertion is immediate, and the second one follows from the fact that the dilation is a unitary operator in $\LL^2(\R^d)$.
\end{proof}

We are now ready to state  the main theorem of this section.

\begin{theorem}\label{wav-set}
The set $\Wr$ is dense in $\LL^2(\R^d)$. 
Precisely, if $f \in \LL^2(\R^d)$, and $\zve >0$, there exists a 
function $\psi \in  \LL^2(\R^d)$, an expansive matrix $a \in \R^{d\times d}$, and a lattice $\Gamma$, such that
\begin{itemize}
\item $\displaystyle{\|f - \psi\|_2 < \zve}$ and
\vspace{1mm}
\item $\displaystyle{(\zC, a ,\zG)}$ is a Riesz basis for $\LL^2(\R^d)$.
\end{itemize}
\end{theorem}

Before proceeding with the proof let us give the following definition.
\bdef
Let $\zG$ be an arbitrary lattice. A set $E \subset \R^d$ is $\zG$-congruent to a set $\tilde E \subset \R^d$, if
there exist partitions $\{E_s: s\in\zG\}$ of $E$ and  $\{\tilde E_s: s\in\zG\}$ of $\tilde E$ in measurable sets such that for every $s\in\zG$, $E_s=\tilde E_s + s.$
\end{definition}

In the proof of Theorem \ref{wav-set} we will use  the $\LL^\infty$-norm of $\R^d.$ In this way we will obtain
an orthogonal basis of exponentials supported on  the cube $B_{\infty}(0,R/2)$ for some appropiate $R$.
If we multiply the elements in this basis by a function that is bounded above and bounded away from zero, by the previous Lemma, we will have a Riesz basis on $\KK_{B_{\infty}(0,R/2)}$. Therefore, looking at the previous construction in the frame case, we
would now need a set $\UU$ on which $|\hat g|$ is bounded away from zero, 
that is $\zG$-congruent to $B_{\infty}(0,R/2)$, and furthermore, that 
tiles the plane by dilations by a matrix $a$. 
This forces us to be more careful in the choice of $r$, and
will also limit our choices of $\zG$ and $a$.

\subsection{Proof of Theorem \ref{wav-set}}
\begin{proof}

Let as before $g \in \LL^2(\R^d)$ be such that $\hat g$ is continuous and
$\|\hat f - \hat g\|_2 < \frac{\zve}{2}$.
We will approximate $\hat g$ with an appropriate function.

\begin{itemize}
\item Let
$R > 0 \text{ be such that} \int_{\R^d \setminus B_{\infty}(0,R/2)} |\hat g(\zo)|^2 d\zo < \frac{\zve^2}{16}$.
\vspace{1mm}
\item
We now select $r>0$ small enough such that:
\begin{align}
 &\quad r < \frac{R}{3} \label{r1} \\
&\quad \int_{B_{\infty}(0,r/2)} |\hat g(\omega)|^2 d\omega  <  \frac{\zve^2}{16} \label{r2} \\
 \text{and} &\quad {\bf \rm m}  r^d <  \frac{\zve^2}{16} \label{r3} \quad 
 \text{where} \quad {\bf \rm m} := \max\{  |\hat g(\zo)|^2: \zo \in B_{\infty}(0,R/2)\}.
\end{align}
\end{itemize}

Let now 
\begin{equation*}
\Gamma = R \Z^d \quad \text{and} \quad a = \frac{R}{r} I_{d\times d}. 
\end{equation*}
We define
\begin{align*}
\widetilde{T}:\R^d &\longrightarrow \R^d\\
   x   & \mapsto x + R \xi_j \quad \text{if}\ x \in j^{\rm th}-\text{quadrant}
\end{align*}
where $\xi_j$ is the vertex of the cube $[-1,1]^d$ that lies in the $j^{\rm th}$-quadrant.

Let us call 
\begin{equation}\label{a0}
A_0 = B_{\infty}(0,r/2) \quad \text{and define}\quad A_i := (a^{-1} \circ \widetilde{T})^i (A_0), i=1, 2,\dots
\end{equation}
It will be convenient to use the notation $A_0^j$ for the intersection of the set $A_0$
with the  $j^{\rm th}$-quadrant. With this notation, note that
\begin{align} A_i & =  
\bigcup_{1\leq j \leq 2^d} (a^{-1}\circ \widetilde{T})^i(A_0^j) \notag\\
                  & =   \bigcup_{1\leq j \leq 2^d}  T_{R a^{-1}\xi_j+\cdots Ra^{-i} \xi_j} (a^{-i}(A_0^j)) . \label{ai}
\end{align}
Here $T_y$ denotes the usual translation by $y$ in $\LL^2(\R^d)$. Therefore we have that
\beq\sum_{i=1}^{\infty} \mu(A_i) = \mu(A_0) \sum_{i=1}^{\infty} \left(\frac{r}{R}\right)^{d i} 
=  r^d \frac{r^d}{R^d-r^d} < r^d,
\label{meas-ai}
\eeq
where the last inequality holds by \eqref{r1}.

We define the set
\beq\notag
\UU  := a\left( \bigcup_{i=0}^{\infty} A_i\right) \setminus \left( \bigcup_{i=0}^{\infty} A_i\right), 
\eeq
and for $\zl = \frac{\zve}{8(R)^{d/2}}$, the function $h$ by
\begin{equation*}
h(\omega) := \begin{cases}
\hat{g}(\omega) & x \in \UU\cap E_\zl \\
\zl & x \in \UU \setminus E_\zl \\
0 & \text{else}
\end{cases}
\end{equation*}
where as before 
\beq \notag
E_{\lambda} :=\{\zo \in \R^d: |\hat g(\zo)| > \lambda \}.
\eeq

Our claim is that $\psi$ with $\hat \psi  = h$ satisfies both conditions of the Theorem.

For the first, let us compute
\beq \notag
\|f -\psi\|_2 = \|\hat{f} - h\|_2 \leq \|\hat{f} - \hat g\|_2 + 
\|\hat g - h\|_2 < \frac{\zve}{2} + \|\hat g-h\|_2.
\eeq
To compute $\|\hat g - h\|_2$ we note that
$\|\hat g - h\|_2^2 $ can be split into 3 integrals
\beq
\int_{\UU\cap E_\zl} |\hat g - h|^2  + \int_{\UU \setminus E_\zl } |\hat g - h|^2 
+ \int_{\R^d\setminus \UU} |\hat g - h|^2 . \label{h-g}
\eeq
By the definition of $h$, the first term in \eqref{h-g} vanishes.

For the second, note that on $\UU\setminus E_\zl$, $|\hat g(\zo)| \leq \zl$ and therefore
\beq \label{h-g2}
 \int_{\UU \setminus E_\zl } |\hat g -  h|^2  \leq 4 \zl^2\mu(\UU\setminus E_{\zl}) < \frac{\zve^2}{16} .
 \eeq
 
 For the last term, since $h(\zo) = 0$ if $\zo \not\in \UU$ we have 
  \beq\label{xx}
  \int_{\R^d\setminus \UU} |\hat g|^2  < \int_{B_{\infty}(0,R/2)\setminus \UU} 
|\hat g|^2  + \int_{\R^d\setminus B_{\infty}(0,R/2)} |\hat g|^2.
  \eeq
The right hand side of (\ref{xx}) can be written,
 \beq
  \int_{A_0} |\hat g|^2  +  \int_{\cup_{i\geq1} A_i} |\hat g|^2  + 
\int_{\R^d\setminus B_{\infty}(0,R/2)} |\hat g|^2 .
  \eeq
The first and last term in this sum are each smaller than $\zve^2/16$ by our choice of $r$ and $R$. For the middle one, we use the computation about the measure of $\cup_{i\geq 1}A_i$ done in \eqref{meas-ai} and the choice of $r$ in \eqref{r3} to obtain
 \beq
 \int_{\cup_{i\geq1} A_i} |\hat g|^2 d\zo \leq {\bf \rm m} \mu(\cup_{i\geq1} A_i) < {\bf \rm m} r^d
 < \frac{\zve^2}{16}.
  \eeq
Putting all this together, we obtain
\beq
\|\hat g - h\|^2_2 < \frac{\zve^2}{4} \quad \text{and therefore} \quad
\|f -\psi\|_2   < \zve.
\eeq

This proves that the function $h$ can be chosen as close to $\hat g$ as we wish. It remains to show, that $(\zp,a,\zG)$ is a a Riesz basis. For this, we observe that:
\begin{itemize}
\item By construction, $\UU$ tiles $\R^d\setminus{0}$ by dilations by $a$.
\item Furthermore $\UU$ tiles $\R^d$ by translations on $\Gamma$. For this, we first note that if $x \in A_n$ then
\beq
\frac{rR}{R-r} \left(1-\left(\frac{r}{R}\right)^n\right) \leq \|x\|_{\infty} \leq 
\frac{rR}{R-r} \left(1-\left(\frac{r}{R}\right)^n\left(\frac{r+R}{2R}\right) \right). \label{size-ai}
\eeq
This fact allows us to conclude that:
\begin{enumerate}
\item $A_i \subset \left(B_{\infty}(0,R/2) \setminus A_0\right), i \geq 1$,
\item $A_i \cap A_j = \emptyset \ i \not= j$
\end{enumerate} 
which allows us to rewrite $\UU$ 
\begin{align*}
\UU & = a\left( \bigcup_{i=0}^{\infty} A_i\right) \setminus \left( \bigcup_{i=0}^{\infty} A_i\right) \\
 & = \left(aA_0 \setminus ( \bigcup_{i=0}^{\infty} A_i)\right) \cup  \left( \bigcup_{i=0}^{\infty} \tilde T A_i\right).
 \end{align*}
This shows that $\UU$ is $\zG$-congruent to $B_{\infty}(0,R/2)$.
\end{itemize}

On the other hand, since
$\zl < |h(\zo)| < {\bf \rm m} $ on $\UU$  by Lemma~\ref{riesz} 
$$ \left\{h(\zo) \frac{1}{R^{d/2}} e^{2\pi i \frac{\bf k}{R}\zo} : k \in \Z^d\right\}$$
is a Riesz basis for $\KK_{\UU}$.

Thus we found a Riesz basis for $\KK_{\UU}$, and $\UU$ is a {\em wavelet set} for the dilation $a$ 
and translation $\zG$. Therefore, the system $(\zp,a,\zG)$ is a Riesz basis of $\LL^2(\R^d)$.

\end{proof}

\begin{remark}
In the above proof, we constructed a function that is supported on a wavelet set. In this way we in fact prove that a proper subset of $\Wr$ is dense in $\LL^2(\R^d)$, since the functions in the dense set are supported on a wavelet set. Further, they generate {\em quasi-orthogonal affine systems}, since they have orthogonal dilations (i.e. 
$<D_{a^j}T_{\zg}\psi, D_{a^{j'}}T_{\zg'}\psi> = 0$ if $j\not= j'$).
\end{remark}

\section{Density and Orthogonal Wavelets}

\label{ort-sys}
In this section we will consider the problem of the density for orthonormal 
wavelets. Since orthonormal wavelet functions have norm one, 
the natural question in this context is wether they are dense on the 
unit sphere of $\LL^2(\R^d)$ with the induced metric 
(i.e. $\Sph_d := \{ f \in \LL^2(\R^d): \|f\|_2 = 1\}$).

An immediate argument, that we will see later, shows that if we fix 
either the dilation matrix $a$ or the lattice $\zG$, then the sets 
$\Wo(a)$ and $\Wo(\zG)$ are not dense in $\Sph_d$. If we allow both, 
the matrix $a$ and the lattice $\zG$ to be arbitrary, it is an open 
problem  if the set $\Wo$ is dense in $\Sph_d$.

In this section we will prove that the set of functions 
\beq \notag
\widetilde{\Wo} := \left\{\psi: \exists \ \zG \text{ and }
a \in \GL, \text{ such that } (\psi, a, \zG) 
\text{ is an  orthonormal system}\right\}
\eeq  
is dense in $\Sph_d$. Note that in $\widetilde{\Wo}$ we removed the completeness requirement of the system.

In particular we prove an interesting property of the Grammian of an 
arbitrary function 
in $\LL^p([0,1]^d)$ which is of independent interest. This property is 
then used to derive some consequences which in particular imply the 
above mentioned result.

Throughout this section, $b$ will be a matrix in \GL.
Our first Lemma is a known result (see for example \cite{Mal89b}).
We state it here in the form we need it. 
\begin{lemma}
Let $f \in \LL^2(\R^d)$.  The following statements are equivalent
\benu
\item \label{i} The system $\{T_{((b^*)^{-1}k)}f : k \in \Z^d\}$ is orthonormal in 
$\LL^2(\R^d)$.
\item \label{i-1} The system $\{\hat f(\zo) e^{-2\pi i \zo\cdot (b^*)^{-1}k}: k \in \Z^d\}$ is orthonormal in 
$\LL^2(\R^d)$.
\item \label{ii} $\sum_{s\in \Z^d} | |\det b|^{1/2} \hat f(b(\zo + s)|^2 = 1$
 a.e. in $[0,1]^d$.
\item \label{iii} $\sum_{s\in \Z^d} | \hat f(\zo + bs)|^2 = |\det b|^{-1}$ 
a.e. in $b([0,1]^d)$.
\eenu
\end{lemma}
\begin{proof}
We will show (\ref{i}) 
$\Leftrightarrow$ (\ref{ii}). The rest is trivial.

The system $\{T_{((b^*)^{-1}k)}f : k \in \Z^d\}$ is orthonormal in $\LL^2(\R^d)$,
if and only if
\begin{align*}
\zd_{0,k} & = \int_{\R^d} f(x-(b^*)^{-1}k) \overline{f(x)} dx =
\int_{\R^d} 
e^{-2\pi i \zo \cdot (b^*)^{-1}k} \hat f(\zo) \overline{\hat f(\zo)} d\zo \\
       &=   \int_{[0,1]^d} {\big (}|\det b|\sum_s |\hat f(b(\zo + s))|^2 {\big )}\ 
       e^{-2\pi i \zo \cdot k}  d\zo. 
\end{align*}
Hence, it will be orthonormal if and only if (\ref{ii}) is satisfied.
\end{proof}

Note that the sum in (\ref{ii}) can be rewritten as  
$\sum_{s\in \Z^d} | T_s D_b \hat f(\zo)|^2$, which motivates the following definition. 
\begin{definition} If $g \in \LL^2(\R^d)$ and $c \in \GL$, the Grammian of $g$ with respect to $c$ is
the function
\beq\notag
g_c(\zo) = {\big (}\sum_{s\in \Z^d} | T_s D_c g(\zo)|^2{\big )}^{1/2} = 
|\det c|^{1/2} {\big (}\sum_{s\in \Z^d} | g(c(\zo + s))|^2 {\big )}^{1/2}
\eeq
Note that $g_c$ is a $\Z^d$-periodic function.
\end{definition}
In what follows the set 
$E(g,b) := \{\zo \in \R^d: g_b(\zo) > 0\}$ will be relevant.
It is immediate to verify that, if
$g \in \LL^2(\R^d)$,
then
\beq \label{norm-gram}
\|g_b\|_{\LL^2([0,1]^d)} = \|g\|_{\LL^2(\R^d)}.
\eeq

For $g \in \LL^2(\R^d)$, we define $u_g$ by
\beq \label{g-normalizada}
u_g(\zo) := \begin{cases}\frac{ g (\zo)}{g_b(b^{-1}\zo)} & \text{ if } \zo \in b E(g, b)\\
0 & \text{otherwise}. 
\end{cases}
\eeq
We have the following Lemma, whose proof is immediate.
\begin{lemma} \label{lem-3}
Let $g \in \LL^2(\R^d)$ and $b\in \R^{d\times d}$ an invertible matrix. 
Then the
function $u_g$ defined in \eqref{g-normalizada} 
satisfies that
\beq \notag 
(u_g)_b(\zo) = |\det b|^{1/2} {\big (}\sum_{s \in \Z^d} |u_g(b(\zo + s))|^2 {\big )}^{1/2} 
= 1\quad \text{a.~e.~} \zo \in  E(g, b).
\eeq
In particular, if $\mu(E(g, b) \cap [0,1]^d) = 1$ then the function $\zvf$ defined by $\hat \zvf = u_g$ has orthonormal\ $(b^*)^{-1}\Z^d$ translates 
(i.e. $\{T_{((b^*)^{-1}k)}\zvf: k \in \Z^d\}$ is orthonormal).
\end{lemma}

The next Lemma states that for each invertible matrix $b$ the distance in 
$\LL^2(\R^d)$ of two arbitrary functions is bigger than the distance of its
Grammians in $\LL^2([0,1]^d)$.

\begin{lemma} \label{lem-4} Let $g, h \in \LL^2(\R^d)$, $b\in \GL$
and ${u_{g}}$ as in \eqref{g-normalizada}. Then we have
\begin{align*}
 \|g - h\|_{\LL^2(\R^d)}   & \geq 
\|{g}_b - {h}_b\|_{\LL^2([0,1]^d)}\quad \text{and}\\
\|g - u_{g}\|_{\LL^2(\R^d)}  & =
\|g_b - \chi_{E(g, b)\cap [0,1]^d} \|_{\LL^2([0,1]^d)}.
\end{align*}

\end{lemma}

\begin{proof}
By \eqref{norm-gram}
\begin{align*}
 \|g-h\|_{\LL^2(\R^d)}^2 &= \|(g - h)_b\|_{\LL^2([0,1]^d)}^2  \\
&= \int_{[0,1]^d} |\det b| \sum_s |g(b(\zo + s)) - h(b(\zo + s))|^2 d\zo \\
&\geq 
\|{g}_b - {h}_b\|_{\LL^2([0,1]^d)},
  \end{align*}
 where in the last inequality we used for the $\ell_2(\Z)$ norm the inequality $\|x -y\| \geq
 |\|x\| - \|y\| |$.

For the second equation, we compute directly
\begin{align*}
\|g - u_{g}\|^2_2 & = 
\int_{b E(g, b)} \left|{g}(\zo) - 
\frac{g(\zo)}{(|\det b|\sum_s |g(\zo + bs)|^2)^{1/2}}\right|^2 d\zo\\
& =    \int_{b E(g,b)} \frac{|g(\zo)|^2}{|\det b|\sum_s |g(\zo + bs)|^2} 
\left|\left(|\det b|\sum_s |g(\zo + bs)|^2\right)^{1/2} - 1\right|^2 d\zo. 
\end{align*}
If we now periodize and change variables, we obtain the result.
\end{proof}

\begin{remark} Note that if $\mu(E(g,b) \cap [0,1]^d) = 1$ then
\beq \notag
\|g- u_{g}\|_{\LL^2(\R^d)}   =
\|g_b - \chi_{[0,1]^d} \|_{\LL^2([0,1]^d)}.
\eeq
\end{remark}

A consequence of Lemma~\ref{lem-4} is the following: Assume that we want to be able to find a function 
$\psi$ that is close to a 
given function $f \in \Sph_d$ and that has orthonormal translates with respect to a lattice $(b^*)^{-1}\Z^d$. Then we will need that $\Gfa$ is close 
(in $\LL^2([0,1]^d)$) to the constant function $1$. Now, if for a given matrix 
$b$, $\Gfa$ is far from $\chi_{[0,1]^d}$, will the choice of a different 
matrix $b$ improve the error? The next theorem establishes the 
interesting result that $\Gfa$ in fact converges almost everywhere and 
also in norm to $\chi_{[0,1]^d}\|f\|_2$ when $\|b\| \rightarrow 0$. 

There is a natural interpretation in  time domain of this result: 
when $\|b\|$ becomes small, then $\|(b^*)^{-1}\|$ grows, that is the associated lattice 
$(b^*)^{-1}\Z^d$ becomes sparser and since our functions are in $\LL^2(\R^d)$ they 
will have some decay at infinity which will imply that the scalar 
product between two of its translates will be small. 

In fact the theorem is more general. The convergence also holds for the $p$-Grammians of a function $f \in \LL^p$ that we
denote by $f_{b,p}$ and are defined as
$ f_{b,p} (\zo) = \left(|\det b|\sum_{s\in \Z^d} | 
f(b(\zo+s))|^p\right)^{1/p}.
$ Note that for $p=2$, $f_{b,2}$  coincides with our previous $f_b$.

The following theorem is a generalization of a result proved in  \cite{JWW05} for the $\LL^1$ 
case in a completely different context.

\begin{theorem}\label{teo-gram}
For any $f \in \LL^p(\R^d)$, 
\beq 
f_{b,p} \longrightarrow \|f\|_p\chi_{[0,1]^d} 
\quad \text{a.e. and in $\LL^p([0,1]^d)$ when}\quad \|b\| \rightarrow 0.
\eeq
\end{theorem}

We postpone the proof of the theorem until the last section.

 Let us now see, how we deduce from Lemma \ref{lem-4} immediately  that when the matrix $b$ is fixed, that is the lattice $\zG = (b^*)^{-1}$ is fixed,
then the set $\Wo(\zG)$ is not dense in the sphere:
By Parseval, Lemma \ref{lem-4} and using that $\hat g_b= 1 $ a.~e.~for every $g\in \Wo(\zG)$ we have for $f \in S_d$ and 
$g\in \Wo(\zG)$:
$$
 \|f - g\|_{\LL^2(\R^d)}  \geq 
\|\hat{f}_b - \hat{g}_b\|_{\LL^2([0,1]^d)} = \|\hat{f}_b -\chi_{[0,1]^d}\|_{\LL^2([0,1]^d)}
$$
The proof is completed choosing a function $f$ in $\Sph_d$ whose Grammian (with respect to $b$) is far from $\chi_{[0,1]^d}$.

The following argument from Yang Wang [private communication] shows that if the dilation $a$ is fixed,
then $\Wo(a)$ is not dense in $\Sph_d$:
Assume that  $\Wo(a)$ is dense in $\Sph_d$. Let $f\in \Sph_d$ be an arbitrary function such that
$ <f,D_af> \not= 0$.
Choose $\psi_n \in \Wo(a)$ such that $\psi_n \rightarrow f, n\rightarrow \infty$ in $\Sph_d$ . Then we have:
$$
0\  =\ <\psi_n,D_a\psi_n>\  \rightarrow\  <f,D_af>,
$$ 
which is a contradiction.

Now we are ready to prove a density result for the set of generators of
orthonormal (not necessarily complete) wavelet systems:
\begin{theorem}\label{teo-osys}
The set $\widetilde{\Wo}$ of generators of orthonormal wavelet systems is dense in $\Sph_d$
in the induced $\LL_2(\R^d)$ metric. 
\end{theorem}

\begin{proof}
Let $f \in \Sph_d$ and $\zve>0$ be given.
Choose $r, R, g$ and $h$ as in Theorem~\ref{teo-2}.

By Theorem~\ref{teo-gram} there exists $b \in GL_d(\R)$ with small enough norm,
such that
\beq \notag
\|h_b - \chi_{[0,1]^d}\|_{\LL_2([0,1]^d)} < \frac{\zve}{2}
\eeq
and $h_b > 0$ a.e. Note that this is possible since $|h| > 0$ in $\{\zo: r\leq \|\zo\|
\leq R\}$. 
Now using Lemma~\ref{lem-4} 
\beq \notag
\|h - u_{h}\|_{\LL^2(\R^d)} < \frac{\zve}{4}
\eeq
 and therefore the inverse
Fourier transform of $u_{h}$
has orthonormal 
translates in the lattice $(b^*)^{-1}\Z^d$.

By choosing the dilation $a$ to be  $a = \frac{R}{r} I_d$, the set
\beq \notag
\left\{
|\det a|^{j/2} u_{h}(a^j \zo) e^{-2\pi a^j \zo (b^*)^{-1} k}: k \in \Z^d, j \in \Z 
\right\}
\eeq
is orthonormal, and consequently, if $\zp$ is such that $\hat \zp = u_{h}$
\beq \notag
\left\{
|\det a|^{j/2} \zp(a^j x - (b^*)^{-1}k): k \in \Z^d, j \in \Z 
\right\}
\eeq
is an orthonormal wavelet system, and so $\zp \in \Wo$.
Further we have
\beq \notag
\|f - \zp\| \leq \|\hat f - \hat g\| + \|\hat g - h\| + \|h - u_h\| < \zve
\eeq
which shows that $\Wo$ is dense in $\LL^2(\R^d)$.
\end{proof}

\subsection{Proof of Theorem~\ref{teo-gram}}
\begin{proof}
We will divide the proof in several steps.
Let us denote the unit cube by $Q = [0,1]^d$. We will denote by $\|\cdot \|_p$ the $p$-norm in $\LL^p(\R^d)$, and
$\| \cdot \|_{\LL^p(Q)}$ the $p$-norm in $\LL^p(Q)$.

\vspace{\baselineskip}
\noindent
$\bullet$
We first prove the theorem for the case that $f = \chi_I$ where $I$ is a finite $d$-dimensional interval. In
this case, we compute
\beq \label{caracteristica}
f_b (\zo) = |\det b| \sum_s \chi_I(b(\zo+s)) =  |\det b| \sum_s \chi_{b^{-1}I} (\zo+s).
\eeq
Let
\begin{align*}
N_i(b) & := \{s \in \Z^d: bQ+bs \subset I\}\\
N_o(b) & := \{s \in \Z^d: \mu((bQ+bs)\cap I) > 0 \text{ and }\mu((bQ+bs)\cap I^c) > 0\},
\end{align*}
and call $n_i = \card N_i$ and $n_o = \card N_o$. Then, 
\beq \label{sum}
\sum_s \chi_{b^{-1}I} (\zo+s) = (n_i + m(b,\zo)),\ m(b,\zo) \in \N, 0\leq m(b,\zo)\leq n_o.
\eeq
Note that $|\det b|\, n_i(b) \leq \mu(I) \leq |\det b|(n_i(b) + n_o(b))|$ and hence
\beq \label{mui}
0\leq \mu(I) - |\det b|\, n_i (b) \leq |\det b|\, n_o(b).
\eeq
We will now see that $|\det b|\,n_o(b) \rightarrow 0$ for $\|b\|\rightarrow 0$.
For this take $\zve >0$, and let $b\in \GL$ be such that $\|b\| < \zve/\diam(Q)$. 
Since $\|b(x-y)\| \leq \|b\|\|x - y \|$ then $\diam (bQ) \leq \|b\|\diam(Q) < \zve$. 
Therefore, if
$s \in N_o(b)$,
$ bQ+bs \subset \{\zo\in \R^d: \dist(\zo, \partial I) < \zve\}$. Now 
\beq
|\det b|\,n_0(b) = \sum_{s\in N_o} \mu(bQ+bs) \leq \mu (\{\zo\in \R^d: \dist(\zo, \partial I) < \zve\}) < \zve c_1.
\eeq
Therefore we see from \eqref{mui} that 
\beq
|\det b|\,n_i(b) \rightarrow \mu(I) \quad \text{when}\quad \|b\|\rightarrow 0,
\eeq
and using \eqref{mui} we have $ |\det b| \sum_s \chi_{b^{-1}I} (\zo+s) \rightarrow \mu (I)$ a.e. Furthermore,
$f_{b,p} \rightarrow \mu(I)^{1/p} = \|f\|_p.$ 

In addition, since  $|\det b| \sum_s \chi_{b^{-1}I} (\zo+s)
\leq |\det b|(n_i(b) + n_o(b)) < C$ for all $b$ such that $\|b\| < \zd <1$ we can apply dominated convergence to obtain 
\beq
\int_Q\left| f_{b,p} - \mu(I)^{1/p} \chi_Q \right|^p \longrightarrow 0 \quad \text{when}\quad \|b\| \rightarrow 0.
\eeq

\vspace{\baselineskip}
\noindent
$\bullet$ Now we will prove the theorem for the case that $f$ is a finite linear combination of characteristic functions of intervals, i.e. $f = \sum_{j=1}^k \za_j \chi_{I_j}$, where $I_j$ are almost disjoint $d$-dimensional intervals (here by almost, we mean that the intersection can have at most measure $0$). By the disjointness of the $I_j$ we immediately obtain
\beq
f_{b,p} (\zo) = \left(\sum_{j=1}^k |\za_j|^p (\chi_{I_j})_{b,p}^p \right)^{1/p}
\eeq 
and by the previous item, this converges to $\|f\|_p$.

In addition, because of the uniform boundedness of $f_{b,p}$ with respect to $b$, we obtain convergence in $\LL^p(Q)$ by the dominated convergence theorem.

\vspace{\baselineskip}
\noindent
$\bullet$ Assume now that $f > 0, f\in \LL^p(\R^d)$. Let $\zve > 0$ and $0 \leq g_n \nearrow f$ a.e., with $g_n$ as in step 2.
\begin{multline}
\left| f_{b,p}(\zo) - \|f\|_p\right| \leq  \\
\left| f_{b,p}(\zo) - (g_n)_{b,p}(\zo) \right| + \left| (g_n)_{b,p}(\zo) -  \|g_n\|_p\right| +
\left|\|g_n\|_p - \|f\|_p\right| \label{juan}
\end{multline}
First note that, as in Lemma~\ref{lem-4} for any $f, g \in \LL^p(\R^d)$ we have $\|f_{b,p} - g_{b,p}\|_{\LL^p(Q)} 
\leq \|f - g\|_{p}$, and therefore
\beq
(g_n)_{b,p} \xrightarrow{\|\cdot \|_{\LL^p(Q)}} f_{b,p}.
\eeq
We can therefore choose a subsequence $g_{n_j}$, such that $(g_{n_j})_{b,p}(\zo) \rightarrow f_{b,p}(\zo)$ a.e. in $Q$ when
$j \rightarrow \infty$. Hence, for a large enough $j$, the first and third term in \eqref{juan} are each less than $\zve/4$.

Now, for a fixed $j$, we can choose $\zd$ such that, if $\|b\| < \zd$, then
$|(g_{n_j})_{b,p}(\zo) - \|g_{nj}\|_p |< \zve/2$ which gives the point-wise convergence.

For the convergence in $\LL^p(\R^d)$, using the previous results, we write
\begin{multline}
\left\| f_{b,p} - \|f\|_p\right\|_{\LL^p(Q)} \leq  \\
\left\| f_{b,p} - (g_n)_{b,p} \right\|_{\LL^p(Q)} + \left\| (g_n)_{b,p} -  \|g_n\|_p\right\|_{\LL^p(Q)} +
\left\|\|g_n\|_p - \|f\|_p\right\|_{\LL^p(Q)}. \label{juan-1}
\end{multline}
By the particular choice of the sequence $g_n$, we have that $g_n \xrightarrow{\|\cdot \|_p} f$. Therefore the first and last term of \eqref{juan-1} go to $0$ when $n\rightarrow \infty$. The middle term goes to $0$ by the previous step.

Since $f_{b,p} = |f|_{b,p}$ the result holds for arbitrary $f \in \LL^p(\R^d)$.
\end{proof}

\subsection{Hilbert Spaces}

The results in this paper carry over to abstract separable Hilbert spaces,
via unitary isomorphisms,
where the translation and dilation operators are replaced by arbitrary 
unitary operators. General systems obtained in this context 
(unitary systems) have been  studied in  detail in \cite{DL98}.

Finally, we would like to mention that M. Bownik \cite{Bow05}, independently of our work and using different techniques has obtained for the case of wavelet frames some similar density results as in section 3.

\section{Acknowledgments}
The research for this paper was carried during our visit at the ESI 
(Erwin Schroedinger Institute) and NuHAG (Numeric and Harmonic Analysis Group at the University of Vienna) 
during the special semester on ``Modern Methods on Harmonic Analysis''. We thank H.~Feichtinger and 
K.-H.~Gr\"ochenig for the invitation and the ESI and NuHAG for the hospitality during our stay. 
We also acknowledge support from an Al$\beta$an-fellowship (EU:E04E031835AR) and Grants
UBACyT X058, X108; PICT 15033, and CONICET 5650.

We thank Yang Wang for fruitful discussions during his visit to the ESI, and  David Larson
for stimulating conversations that originated this research.

We also thank the anonymous referees for constructive remarks and suggestions that helped to improve the final version of this manuscript.


\newcommand{\etalchar}[1]{$^{#1}$}
\providecommand{\bysame}{\leavevmode\hbox to3em{\hrulefill}\thinspace}
\providecommand{\MR}{\relax\ifhmode\unskip\space\fi MR }
\providecommand{\MRhref}[2]{%
  \href{http://www.ams.org/mathscinet-getitem?mr=#1}{#2}
}
\providecommand{\href}[2]{#2}


\end{document}